\documentclass{gtart_h}  

\def\ifplaintex{\expandafter\ifx\csname documentclass\endcsname\relax}

\def\ifplaintex{\expandafter\ifx\csname documentclass\endcsname\relax}

\ifplaintex 
\else
\fi

\expandafter\ifx\csname epsfbox\endcsname\relax\input epsf\fi

\def\gt{{\mathsurround=0pt\it $\cal G\mskip-2mu$eometry \&\ 
$\cal T\!\!$opology}}        

\def\gtp{{\mathsurround=0pt\it $\cal G\mskip-2mu$eometry \&\ 
$\cal T\!\!$opology $\cal P\!$ublications}}  

\def\lognumber#1{\def\thelognumber{#1}}
\def\volumenumber#1{\def\thevolumenumber{#1}}
\def\papernumber#1{\def\thepapernumber{#1}}
\def\volumeyear#1{\def\thevolumeyear{#1}}

\def\pagenumbers#1#2{\def\startpage{#1}\def\finishpage{#2}}
\def\published#1{\def\publishdate{#1}}
\def\proposed#1{\def\theproposer{#1}}
\def\seconded#1{\def\theseconders{#1}}
\def\received#1{\def\receiveddate{#1}}
\def\revised#1{\def\reviseddate{#1}}
\def\accepted#1{\def\accepteddate{#1}}

\def\coverauthors#1{\def\thecoverauthors{#1}}
\def\asciiauthors#1{\def\theasciiauthors{#1}}
\def\asciiaddress#1{\def\theasciiaddress{#1}}
\def\asciiemail#1{\def\theasciiemail{#1}}

\long\def\asciiabstract#1{\long\def\theasciiabstract{#1}}
\def\asciikeywords#1{\def\theasciikeywords{#1}}
\def\shortauthors#1{\def\theshortauthors{#1}}

\let\\\par\let\thelognumber\relax
\let\thevolumenumber\relax\let\thepapernumber\relax
\let\thevolumeyear\relax\let\thesamplenumber\relax\let\startpage\relax
\let\finishpage\relax\let\publishdate\relax\let\receiveddate\relax
\let\reviseddate\relax\let\accepteddate\relax\let\theasciititle\relax
\let\theasciiauthors\relax\let\theasciiaddress\relax
\let\theasciiabstract\relax\let\theasciikeywords\relax
\let\theasciiemail\relax\let\theshortauthors\relax\let\theshorttitle\relax
\let\thecoverauthors\relax

\long\def\maketitlep{   

\count0=\startpage

\gt\hfill      
\hbox to 77pt{\vbox to 0pt{\vglue -15pt\epsfbox{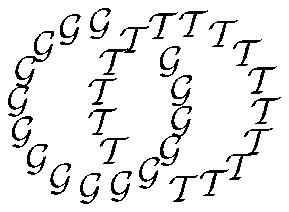}\vss}\hss}
\break
{\small\ifx\thesamplenumber\relax 
Volume \else Sample
\fi\thevolumenumber\ (\thevolumeyear)
\startpage--\finishpage\nl
Published: \publishdate}
\vglue 0.5truein plus 0.4fil minus 0.1truein

{\parskip=0pt\leftskip 0pt plus 1fil\def\\{\par\smallskip}{\ifplaintex\large
\else\Large\fi\bf\thetitle}\par\medskip}   

\vglue 0pt plus 0.1fil 

{\parskip=0pt\leftskip 0pt plus 1fil\def\\{\par}{\sc\theauthors}
\par\medskip}

\vglue 0pt plus 0.1fil 

{\small\parskip=0pt\let\newline\\
{\leftskip 0pt plus 1fil\def\\{\par}{\sl\theaddress}\par}
\expandafter\ifx\theemail\relax    
\relax\else\vglue 5pt plus 0.02fil minus 2pt\def\\{\stdspace{\rm 
and}\stdspace} 
\cl{Email:\stdspace\tt\theemail}\fi
\ifx\theurl\relax                  
\relax\else\vglue 5pt plus 0.02fil minus 2pt\def\\{\stdspace{\rm 
and}\stdspace}
\cl{URL:\stdspace\tt\theurl}\fi\par}

\vglue 7pt plus 0.3fil minus 3pt

{\bf Abstract}
\vglue 5pt plus 0.1fil minus 2pt

\theabstract

\vglue 7pt plus 0.3fil minus 3pt

{\bf AMS Classification numbers}\quad Primary:\quad \theprimaryclass

Secondary:\quad \thesecondaryclass

\vglue 5pt plus 0.3fil minus 2pt

{\bf Keywords:}\quad \thekeywords

\vglue 10pt plus 0.5fil minus 5pt

{\small  Proposed: \theproposer\hfill Received: \receiveddate\nl
Seconded: \theseconders\hfill 
\ifx\reviseddate\relax                         
Accepted: \accepteddate                        
\else
Revised: \reviseddate                          
\fi}
\eject
}       

\font\phead=cmsl9 scaled 950
\font=cmsl9 scaled 1050
\font\pnum=cmbx10 scaled 913
\font\lnum=cmbx10 
\font\pfoot=cmsl9 scaled 950
\font=cmsl9 scaled 1050
\ifplaintex
\headline{\vbox to 0pt{\vskip -4.5mm\line{\small\phead\ifnum
\count0=\startpage ISSN 1364-0380 (on line)
1465-3060 (printed) \hfill {\pnum\folio}\else\ifodd\count0\def\\{ }%
\ifx\theshorttitle\relax\thetitle\else\theshorttitle\fi\hfill{\pnum\folio}
\else\def\\{ and }{\pnum\folio}\hfill\ifx\theshortauthors\relax\theauthors
\else\theshortauthors\fi\fi\fi}\vss}}
\footline{\vbox to 0pt{\vglue 0mm\line{\small\pfoot\ifnum\count0=\startpage
\copyright\ \gtp\hfill\else
\gt, Volume \thevolumenumber\ (\thevolumeyear)\hfill\fi}\vss
}}
\else
\makeatletter
\def\@oddhead{{\small\lhead\ifnum\count0=\startpage ISSN 1364-0380 (on line)
1465-3060 (printed) \hfill {\lnum\number\count0}\else\ifodd\count0
\def\\{ }\ifx\theshorttitle\relax \thetitle \else\theshorttitle\fi\hfill
{\lnum\number\count0}\else\def\\{ and }{\lnum\number\count0}
\hfill\ifx\theshortauthors\relax 
\theauthors\else\theshortauthors\fi\fi\fi}}\def\@evenhead{\@oddhead}
\def\@oddfoot{\small\lfoot\ifnum\count0=\startpage\copyright\ \gtp\hfill\else
\gt, Volume \thevolumenumber\ (\thevolumeyear)\hfill\fi}
\def\@evenfoot{\@oddfoot}
\makeatother
\fi

\newwrite\gtoutfile
\long\gdef\makeheadfile{  
{\def\\{, }\def\s{ }
\immediate\openout\gtoutfile head.xxx
\immediate\write\gtoutfile{Proxy-for: \ifx\theasciiauthors\relax
\theauthors\else\theasciiauthors\fi\s<\ifx\theasciiemail\relax\theemail\else\theasciiemail\fi>}
\immediate\write\gtoutfile{\noexpand\\}
\immediate\write\gtoutfile{Authors: \ifx\theasciiauthors\relax
\theauthors\else\theasciiauthors\fi}
{\def\\{ }\immediate\write\gtoutfile{Title: \ifx\theasciititle\relax
\thetitle\else\theasciititle\fi}}
\immediate\write\gtoutfile{Subj-class: GT or SG or MG etc}
\immediate\write\gtoutfile{MSC-class: \theprimaryclass\ifx\thesecondaryclass\relax\else, \thesecondaryclass\fi}
\immediate\write\gtoutfile{Journal-ref: Geom. Topol. \thevolumenumber
(\thevolumeyear) \startpage-\finishpage}
\immediate\write\gtoutfile{Comments: Published by Geometry and Topology at}
\immediate\write\gtoutfile{\s\s http://www.maths.warwick.ac.uk/gt/GTVol\thevolumenumber/paper\thepapernumber.abs.html}
\immediate\write\gtoutfile{\noexpand\\}
\immediate\write\gtoutfile{}
\ifx\theasciiabstract\relax
\immediate\write\gtoutfile{\theabstract}\else
\immediate\write\gtoutfile{\theasciiabstract}\fi
\immediate\write\gtoutfile{}
\immediate\write\gtoutfile{\noexpand\\}
\immediate\write\gtoutfile{}
\immediate\closeout\gtoutfile}}  

\def\maketitlepage{\maketitlep\makeheadfile}
\let\maketitle\maketitlepage

\lognumber{616}
\received{25 May 2005}
\volumenumber{9}\papernumber{53}\volumeyear{2005}
\pagenumbers{2303}{2317}   
\revised{2 December 2005}
\accepted{3 December 2005}
\proposed{Robion Kirby}
\seconded{Peter Teichner, Cameron Gordon}
\published{10 December 2005}

\usepackage{amsmath,amssymb,graphicx,labelfig}

\def\figref#1{\hyperlink{#1anchor}{Figure~\ref*{#1}}}
\def\anchor#1{\noindent\hypertarget{#1anchor}{\smash{$\phantom{99}$}}\newline}

\newtheorem{thm}{Theorem}[section]
\newtheorem{lemma}[thm]{Lemma}

\theoremstyle{definition}
\newtheorem{conj}[thm]{Conjecture}

\newtheorem{rem}[thm]{Remark}

\newcommand{\C}{{\mathbb C}}

\newcommand{\MH}{{\mathcal H}}
\newcommand{\M}{{\mathcal M}}
\newcommand{\Z}{{\mathbb Z}}

\def\la{\longrightarrow}

\def\q{\quad}
\def\us{\underset}
\def\s{\sigma}
\def\SS{\Sigma}

\def\tn{\textnormal}

\def\bd{\partial}

\begin{document}

\title{Universal manifold pairings and positivity}

\authors{Michael H Freedman$^1$, Alexei Kitaev$^2$, Chetan
  Nayak$^{1,3}$\\Johannes K Slingerland$^1$, Kevin Walker$^1$ and Zhenghan Wang$^4$}
\asciiauthors{Michael H Freedman, Alexei Kitaev, Chetan
  Nayak, Johannes K Slingerland, Kevin Walker and Zhenghan Wang}
\coverauthors{Michael H Freedman, Alexei Kitaev, Chetan
  Nayak\\Johannes K Slingerland, Kevin Walker and Zhenghan Wang}

\shortauthors{Freedman, Kitaev, Nayak, Slingerland, Walker and Wang}

\address{{\rm 1:}\qua Microsoft Research, 1 Microsoft Way, Redmond, WA 98052, USA
\\
{\rm 2:}\qua California Institute of Technology, Pasadena, CA 91125, USA
\\
{\rm 3:}\qua Department of Physics and Astronomy, UCLA, CA 90095-1547, USA
\\
{\rm 4:}\qua Dept of Mathematics, Indiana University, Bloomington, IN
47405, USA\\\smallskip\\
{\tt\mailto{michaelf@microsoft.com}, \mailto{kitaev@iqi.caltech.edu}, \mailto{nayak@physics.ucla.edu}\\\mailto{joost@microsoft.com}, \mailto{kwalker@microsoft.com}, \mailto{zhewang@indiana.edu}}}

\asciiemail{michaelf@microsoft.com, kitaev@iqi.caltech.edu, nayak@physics.ucla.edu, joost@microsoft.com, kwalker@microsoft.com, zhewang@indiana.edu}

\asciiaddress{MHF,CN,JKS,KW: Microsoft Research, 1 Microsoft Way, Redmond, WA 98052, USA
\\
AK: California Institute of Technology, Pasadena, CA 91125, USA
\\
CN: Department of Physics and Astronomy, UCLA, CA 90095-1547, USA
\\
ZW: Dept of Mathematics, Indiana University, Bloomington, IN}

\begin{abstract}
Gluing two manifolds $M_1$ and $M_2$ with a common boundary $S$ yields
a closed manifold $M$.  Extending to formal linear combinations $x=\SS
a_i M_i$ yields a sesquilinear pairing $p=\langle\,\, , \, \rangle$
with values in (formal linear combinations of) closed manifolds.
Topological quantum field theory (TQFT) represents this universal
pairing $p$ onto a finite dimensional quotient pairing $q$ with values
in $\C$ which in physically motivated cases is positive definite. To
see if such a ``unitary" TQFT can potentially detect any nontrivial
$x$, we ask if $\langle x, x \rangle \neq 0$ whenever $x \neq 0$.  If
this is the case, we call the pairing $p$ positive. The question
arises for each dimension $d=0,1,2, \ldots$. We find $p(d)$ positive
for $d=0,1,$ and $2$ and not positive for $d=4$.  We conjecture that
$p(3)$ is also positive. Similar questions may be phrased for
(manifold, submanifold) pairs and manifolds with other additional
structure.  The results in dimension 4 imply that unitary TQFTs cannot
distinguish homotopy equivalent simply connected 4--manifolds, nor can
they distinguish smoothly $s$--cobordant 4--manifolds.  This may
illuminate the difficulties that have been met by several authors in
their attempts to formulate unitary TQFTs for $d=3+1$.
There is a further physical implication of this paper.  Whereas
3--dimensional Chern--Simons theory appears to be well-encoded within
2--dimensional quantum physics, e.g. in the fractional quantum Hall
effect, Donaldson--Seiberg--Witten theory cannot be captured by a
3--dimensional quantum system.  The positivity of the physical Hilbert
spaces means they cannot see null vectors of the universal pairing;
such vectors must map to zero.
\end{abstract}

\asciiabstract{Gluing two manifolds M_1 and M_2 with a common boundary
S yields a closed manifold M.  Extending to formal linear combinations
x=Sum_i(a_i M_i) yields a sesquilinear pairing p=<,> with values in
(formal linear combinations of) closed manifolds.  Topological quantum
field theory (TQFT) represents this universal pairing p onto a finite
dimensional quotient pairing q with values in C which in physically
motivated cases is positive definite. To see if such a "unitary" TQFT
can potentially detect any nontrivial x, we ask if <x,x> is non-zero
whenever x is non-zero.  If this is the case, we call the pairing p
positive. The question arises for each dimension d=0,1,2,.... We find
p(d) positive for d=0,1, and 2 and not positive for d=4.  We
conjecture that p(3) is also positive. Similar questions may be
phrased for (manifold, submanifold) pairs and manifolds with other
additional structure.  The results in dimension 4 imply that unitary
TQFTs cannot distinguish homotopy equivalent simply connected
4-manifolds, nor can they distinguish smoothly s-cobordant
4-manifolds.  This may illuminate the difficulties that have been met
by several authors in their attempts to formulate unitary TQFTs for
d=3+1.  There is a further physical implication of this paper.
Whereas 3-dimensional Chern-Simons theory appears to be well-encoded
within 2-dimensional quantum physics, eg in the fractional quantum
Hall effect, Donaldson-Seiberg-Witten theory cannot be captured by a
3-dimensional quantum system.  The positivity of the physical Hilbert
spaces means they cannot see null vectors of the universal pairing;
such vectors must map to zero.}


\asciikeywords{Manifold pairing, unitary, positivity, TQFT, s-cobordism}

{\small\maketitle}

\section{Introduction}
We begin by establishing notation.  We will work with oriented,
compact, possibly disconnected, smooth manifolds, although some
comments will also be made concerning the unoriented case.  The
choice of smooth category might be essential:  Our vector $x$ is
constructed from a counterexample to the $s$--cobordism theorem,
which is still open in the topological category.

Let $S$ be a $d-1$ dimensional manifold and let $\M_S$ be the $\C
$--vector space of (finite) formal combinations of manifolds $M_i$
with $\bd M_i =S$, so $x=\underset{i}{\SS} a_i M_i  \in \M_S$.
(Note:  If $S$ does not bound, dimension $(\M_S) =0$.)  If we
denote $S$ with the opposite orientation by $\overline{S}$, then
we have a bilinear pairing
\begin{equation}\label{0.1}
    \M_S \times \M_{\overline{S}} \la \M
\end{equation}
given  by $(\us{i}{\SS} a_i M_i, \us{j}{\SS} b_j N_j) \la
\us{i,j}{\SS} a_i b_j M_i \cup_S N_j$, where $\M = \M_\emptyset$
is the vector space of formal linear combinations of closed
$d$--manifolds. To fit better with the role of Hilbert space in
physics we choose to rewrite (\ref{0.1}) as a sesquilinear pairing
\begin{equation}\label{0.2}
     \langle \, \, ,\, \rangle \co \M_S \times \M_{S} \la \M, 
     \,\,\langle \us{i}{\SS} a_i M_i, \us{j}{\SS} b_j N_j \rangle =
     \us{i,j}{\SS} a_i \overline
     {b}_j M_i \cup_S \overline{N}_j
\end{equation}
which is linear in the first entry and conjugate linear in the
second. The map from $\M_S\times\M_{\overline{S}}$ to
$\M_S\times\M_S$ which intertwines between the pairings is just
the conjugate linear extension of orientation reversal on the
second factor.

\begin{figure}[ht!]\anchor{arcfig}
\cl{\includegraphics[width=8cm, height=2cm]{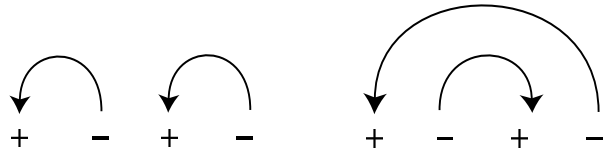}}
 \caption{\label{arcfig} Two inequivalent examples of $1$--manifolds
 with boundary.  In both cases, the manifold consists of two oriented
 line segments and the boundary of two positively and two negatively
 oriented points, but the attaching maps are different}
\end{figure}

We need to be perfectly clear about when two boundary manifolds
$M_i$ and $M_j$ are considered the same element of $\M_S$.  A
basis element $M_i$ of $\M_S$ is a manifold $M_i$ together with a
diffeomorphism  $f_i$ of $\bd M_i$ to $S$.  We say $(M_i ,
f_i)$ and $(M_j , f_j)$  are equivalent if there is a
diffeomorphism \mbox{$\phi\co  M_i \la M_j$} such that:
\begin{equation}\label{0.3}
    f_j \circ \phi \mid_{\bd M_{i}} = f_i.
\end{equation}
With this definition, we have examples where the manifolds $M_i,
M_j$ are the same, but attached differently to the boundary and
hence not equivalent. Perhaps the simplest of these is shown in
\figref{arcfig}, where the manifolds both consist of two line
segments, attached to the four boundary points in different ways.
Less trivially, a surface bounds infinitely many distinct handle
bodies parameterized by the cosets: $MC_g / HC_g$; the genus $=g$
mapping class group modulo the subgroup which extends over a fixed
handlebody.

Occasionally we consider simply the set of bounded manifolds up to
equivalence (ie, the basis vectors of $\M_S$) and denote this set
by $\dot{\M}_S$. We reserve the dot to mean ``unlinearize''.

Our definitions easily extend to (manifold, submanifold) pairs (if
$K\subset M$ is a submanifold we always assume $\bd K \subset
\bd M$). Let $\M^{d, k}$ be the space of formal combinations of
($d$--manifold, $k$--submanifold) pairs.  If $(S, L)$ is a fixed
($(d-1)$--manifold, $(k-1)$--submanifold) pair, we may define
$\M^{d, k}_{(S, L)}$ to be formal combinations of bounding
($d$--manifold, $k$--submanifold) pairs with an equivalence relation
analogous to (\ref{0.3}) and a sesquilinear pairing:
\begin{equation}
    \M^{d,k}_{S, L} \times \M^{d,k}_{S,L} \la \M^{d,k}
\end{equation}
by a formula like (\ref{0.2}).

A variant on gluing pairs is to require the outer manifolds to be
as simple as possible, spheres and disks. This gives sesquilinear
``tangle pairings":
\begin{equation}
    \Im^{d,k}_{L} \times \Im^{d,k}_{L} \la \mathcal{L}^{d,k}
\end{equation}
where $L$ is a fixed $(k-1)$--submanifold of $S^{d-1}$ and $\Im^{d,
k}_L$ is the span of $k$--submanifolds in $D^d$ bounded by $L$. The
target $\mathcal{L}^{d, k}$ is the span of $k$--submanifolds in
$S^d$.

For all the sesquilinear pairings above we may ask if they are
positive, that is, whether $\langle x,x \rangle =0$ implies $x=0$.
The motivation is to understand how much of manifold topology can
potentially be detected by unitary topological quantum field
theories (UTQFTs. See \cite{A} for a definition). To touch on only
the most elementary aspect of this structure, a UTQFT should
assign a scalar to a closed $d$--manifold and a finite dimensional
Hilbert space $V_{S}$ to each $(d-1)$--manifold $S$.  For $X$ with
$\bd X=S$, a vector $\widetilde{X} \in V_{S}$ is assigned and
if $X'$ also satisfies $\bd X' =S$ then $\langle \widetilde{X},
\widetilde{X}'\rangle$ must equal the closed manifold invariant of
$X \overline{X}':=X \cup_S \overline{X}'$. Clearly if one of our
pairings is not formally positive, there will be an $x=\SS a_i
X_i\neq 0$ for which $\langle x,x \rangle =0$, and no unitary TQFT
will be able to distinguish the combination $x$ from zero. This
question is (roughly) in the same spirit as asking if the Jones
polynomial detects all knots.

To make the connection to TQFTs more exact one might choose to
enhance our manifolds with framings, spin structure,
$p_1$--structures, etc$\ldots$, the necessary input for certain
TQFTs.  But the investigation is at such a preliminary stage that
this level of detail is not yet warranted. Also, one may note that
the invariants for closed manifolds often depend only weakly on
the extra structures in the definition of TQFTs, so our results
for closed manifolds may already be useful in such cases. With the
definition complete, let us do what is easy.

\section{Lowest dimensions}

\begin{thm} \label{+pairings}
The following pairings are positive:
\[
\M^d_Y \times \M^d_Y \la \M^d,
\]
for $d=0,1,2$,
\[
 \M^{d,k}_{Y, L} \times \M^{d,k}_{Y,L} \la \M^{d,k},
\]
for $d=0,1,2$; $\, k<d$ and
\[
 \Im^{d,k}_{L} \times \Im^{d,k}_{L} \la \mathcal{L}^{d,k},
\]
for $d=0,1,2$; $\, k<d$.
\end{thm}

\begin{conj} Theorem \ref{+pairings} extends to $d=3$ in the above
cases.
\end{conj}

\begin{lemma} \label{cxlemma}
Suppose there exists a function (the ``complexity function")\break $C\co 
\dot{\M}^d \to \cal{V}$, where $\cal{V}$ is some partially ordered
set, such that for all $M, N \in \dot{\M}^d_Y$, $M \neq N$ implies
$C(M \cup \overline{N}) < 
\max(C(M \cup \overline{M}), C(N \cup \overline{N}))$. 
Then the pairing for
$\M^d_Y$ is positive. Similar statements hold for
$\M^{d,k}_{Y,L}$ and $\Im^{d,k}_{L}$.
\end{lemma}

\proof The hypothesis of the lemma implies that
the terms of maximal complixity in the right-hand side of Equation
\ref{0.2} all lie on the diagonal. Since all coefficients on the
diagonal are positive, there can be no cancellation among these
terms. \endproof

\proof[Proof of Theorem \ref{+pairings}] By the previous
Lemma, it suffices in each case to define an appropriate
complexity function $C$.

We ignore $\mathcal{L}^{d,k}$, since these cases are implied by
the $\M^{d,k}$ cases.

For $\M^1$, define $C(M)$ to be the number of components of $M$.
Let $Y$ be the 0--manifold with $j$ positive and $j$ negative
points. Let $M, N \in \M^1_Y$, and assume for the moment that
neither $M$ nor $N$ contain closed components. Then $C(M \cup
\overline{M}) = C(N \cup \overline{N}) = j$. If $M \neq N$ then at
least one component of $M \cup \overline{N}$ contains 4 or more
arcs, and so $C(M \cup \overline{N}) < j$. Then general case
(where $M$ and $N$ might have closed components) is similar.

Next consider $\M^2_Y \times \M^2_Y \la \M^2$, where $Y$ is the
disjoint union of $j$ circles. Let $M \in \dot{\M}^2$ (ie, $M$ is
a closed oriented 2--manifold). Let $n$ be the number of connected
components of $M$, let $\chi$ be the Euler characteristic of $M$,
and let $\chi_1, \ldots, \chi_n$ be the Euler characteristics of
the components of $M$, listed in increasing order. Define the
complexity of $M$ to be the lexicographic tuple
\[
    C(M) = (n, -\chi, -\chi_1, \ldots, -\chi_n).
\]
The smallest integer that can appear in the tuple is $-2$ so we
formally pad tuples by adding a list of $-3$'s at the end so that
tuples of different lengths can be compared.

Now let $M, N \in \dot{\M}^2_Y$ and assume $M \neq N$. For
simplicity, assume that neither $M$ nor $N$ contain any closed
components. $M$ determines a partition of the components of $Y$:
two components of $Y$ are in the same part of the partition if
they and connected by a component of $M$. In the same way $N$ also
determines a partition of the components of $Y$. If these two
partions differ then $M \cup \overline{N}$ has fewer components
than at least one of $M \cup \overline{M}$ and $N \cup
\overline{N}$, so the hypothesis of the lemma is satisfied. We
assume from now on that the partitions associated to $M$ and $N$
are the same. If $\chi(M) < \chi(N)$, then $C(M \cup \overline{N})
< C(M \cup \overline{M})$ (consider the second component of the
complexity tuple), so we assume from now on that $\chi(M) =
\chi(N)$. The components of $M$ and $N$ are paired according to
their common partition of the components of $Y$. Since $M \neq N$,
there must be at least one pair with differing Euler
characteristics. Amongst such components in non-matching pairs,
choose the one with lowest Euler characteristic, and assume WLOG
that this extremal component belongs to $M$. It follow that $C(M
\cup \overline{N}) < C(M \cup \overline{M})$, and we are done.

The remaining cases of Theorem \ref{+pairings} are proved
similarly. The most complicated case is $\M^{2,1}$. For $(M, K)
\in \dot{\M}^{2,1}$ define $C(M, K)$ to be the lexicographic
triple $(C(M), C(K), C(M \setminus K))$, where $C$ is as above for
plain 1-- and 2--manifolds. \endproof

\begin{rem} Because the Turaev--Viro TQFTs do not
require orientations it is reasonable to also investigate the
universal pairings in the context of unoriented manifolds.  In
this context, define $\langle \us{i}{\SS} a_i M_i, \us{j}{\SS}b_j
N_j \rangle = \us{i,j}{\SS} a_i \overline{b}_j M_i N_j \rangle$.
Theorem \ref{+pairings} continues to hold in the unoriented
context.
\end{rem}

\section{Three dimensional pairings}

This might be the most interesting case and we hope it will be the
subject of future research.  We establish positivity only in a few
rather easy cases where all the work is contained in old theorems.

\begin{thm}\label{+knots}
The pairings $\M^3_S \times \M^3_S \la \M^3$, $S$ a fixed
(possibly empty) union of $\,2$--spheres, and \mbox{$\Im^{3,1}_{2}
\times \Im^{3,1}_{2} \la \mathcal{L}^{3,1}$}, where $2$ denotes
two points with opposite orientations, are positive.
\end{thm}

\proof The essential ingredient in both arguments
is the existence and uniqueness of prime decompositions of
$3$--manifolds \cite{M} and knots \cite{S}. Using this, both cases
reduce to the following lemma:
\begin{lemma} \label{ringlem}
Consider the polynomial rings $\C[p_i \ldots p_n]$ on
indeterminates\break $p_1,\ldots,p_n$ and with a fixed antilinear
involution $^-$ which sends indeterminates to indeterminates. The
natural sesquilinear pairing on these rings,
\begin{equation*}
\begin{array}{rclcl}
\C[p_1 \ldots p_n] &\otimes & \C [p_1 \ldots p_n] &\la& \C[p_1 \ldots p_n]\\
a   &  \otimes & b &\mapsto& \q a\overline{b},
\end{array}
\end{equation*}
is formally positive.
\end{lemma}

\proof We define the complexity of a monomial as
some kind of list of prime powers it contains. If two distinct
primes are related by involution we form a lexicographic pair:
(sum of the two exponents, the smaller of the two exponents).  For
primes paired with themselves the pair is simply \mbox{(exponent,
zero)}. Note that the latter case will arise in the proof of
theorem \ref{+knots}, since some prime knots ($3$--manifolds) are
diffeomorphic to their arrow reversed (orientation reversed)
mirror image.

Now list the pairs in order (padded by $(0,0)$'s) and use this
list of pairs to lexicographically order monomials. Suppose
$X_i=\prod_{k} p^{d_{i, k}}_{i,k}$ are monomials, $x=\SS a_i X_i$,
 then the monomials of greatest complexity in $\langle x,x\rangle$
 are among the diagonal terms $a_i \overline{a}_i X_i \overline{X}_i = a_i
 \overline{a}_i \prod_k p^{d_{i, k}}_{i,k} \overline
 {p}^{d_{i,k}}_{i,k}$. This is easily checked. \endproof

\noindent We make one further observation:
\begin{thm}\label{3.2}
Consider the pairing $\M^3_S \otimes \M^3_S \la \M^3$ where $S$ is
a connected $\tn{genus} = g$ surface. Suppose $x$ is a linear
combination $\SS a_i X_i$, where each $X_i$ is a handle body of
$\tn{genus} = g$ (but each attached to $S$ differently), then
$\langle x,x \rangle =0$ implies each $a_i =0$.
 \end{thm}

\proof The diagonal terms yield $\us{i}{\SS} a_i
\overline{a}_i \left(S^1 \times S^2 \# \cdots \# S^1 \times S^2
\right)$, a positive multiple of the double of the $\tn{genus} =
g$ handle body.  Could an off diagonal term cancel this
contribution? Such an off diagonal term must have $X_i
\overline{X}_j$ diffeomorphic to $N_g := \us{g~\tn{copies}}{\#}
S^1 \times S^2$, and in fact constitute an {\em exotic} minimal
Heegaard decomposition of $N_g$ (if not exotic, our equivalence
relation on bounding manifolds (\ref{0.3}) --- allows us to write
the $X_i \overline{X}_j$ term as an untwisted double). This
contradicts a theorem of Waldhausen's \cite{W} which shows that
any two Heegaard decompositions of $N_g$ are related by an
(orientation preserving) diffeomorphism.  \endproof

\begin{rem} $3$--manifold topology, in practice, is often two subjects,
``comp\-ress\-ion-theory" and hyperbolic geometry, patched together.
Some geometric arguments (about volumes, lengths, and Ricci flow \cite{AST})
offer hope for positivity on the hyperbolic side of the subject
and Theorem \ref{3.2} offers hope on the compression side. For
these reasons we conjecture that the three dimensional pairings
are positive.\end{rem}

\section{4--manifold pairings}
For a variety of $3$--manifolds $S$, we can exhibit vectors $x=M-M'
\in \M^4_S$ such that
\begin{equation}
\langle x, x\rangle = M \overline{M} - M \overline{M}' -
M'\overline{M} + M' \overline{M}' = 0 \in \M^4 =\M^4_\emptyset.
\end{equation}
In all cases the difference between $M$ and $M'$ is a matter of
differentiable structure on an underlying Poincar\'{e} pair (or,
when $\pi_1(M)=\{e\}$, an underlying topological manifold). In one
example of such an $x$, which has its roots in \cite{Ak,AkK1}, $M$
and $M'$ are both copies of the ``Mazur manifold" and $S=\bd
M$, but $M$ is attached by the ``identity" and $M'$ is attached by
a diffeomorphism $\theta$ of the boundary which does not extend to a
diffeomorphism of the interior (but does extend as a
homeomorphism).  According to the definition in equation
(\ref{0.3}), $M$ and $M'$ are distinct, so $x\neq 0$.

In \cite{AkK1} Akbulut and Kirby showed, by direct handle manipulation, that
doubling Mazur's contractible manifold (via the identity on its
integral homology sphere boundary) yields the smooth $4$--sphere,
$M\overline{M}\cong S^{4}$ and remarkably, the $\theta$--twisted
double is also diffeomorphic to the $4$--sphere,
$M\overline{M}'=S^4$ (since $\theta^2 = \mathrm{id}$, $M'\overline{M}\cong
M\overline{M}'$ and trivially $M'\overline{M}'\cong M\overline{M}$). In
\cite{F}, one of us showed that $M$ and $M'$ constituted the same
topological manifold structure on the Poincar\'{e} pair. Then,
with the introduction of gauge theory in topology, it became
possible to distinguish $M$ and $M'$ as smooth structures. Akbulut
\cite{Ak} did this by showing that the Kummer surface $K$ and one
of its logarithmic transforms $K'$, although distinguished by
Donaldson invariants, differed on a combinatorial level only by
cutting out an embedded $M$ and regluing it via $\theta$. This
shows that $M\neq M'\in \M^{4}_{S}$. In fact, if $C$ is the closed
component $C=\overline{K\setminus M}$, we may write
$$K=C\cup M,~~~~K'=C\cup M'$$
$$K\not\cong K' \Longrightarrow M \neq M'.\leqno{\hbox{so}}$$
The pair of manifolds $K,K'$ is but one of the many examples of
pairs of (smoothly) $h$--cobordant but non-diffeomorphic manifolds.

\begin{figure}[ht!]\anchor{3condfig}
\cl{
\small
\SetLabels
\E(.5*.95) $M'$\\
\E(.5*.05) $M$\\
\E(1*.86) $Q$\\
\E(1*.5) $W$\\
\E(1*.15) $P$\\
\E(.31*.3) $N$\\
\E(.43*.8) $2$\\
\E(.43*.58) $3$\\
\E(.6*.2) $2$\\
\E(.6*.41) $3$\\
\L(.52*.54) $X\cong B^5$\\
\endSetLabels
\AffixLabels{
\includegraphics[width=9cm, height=5cm]{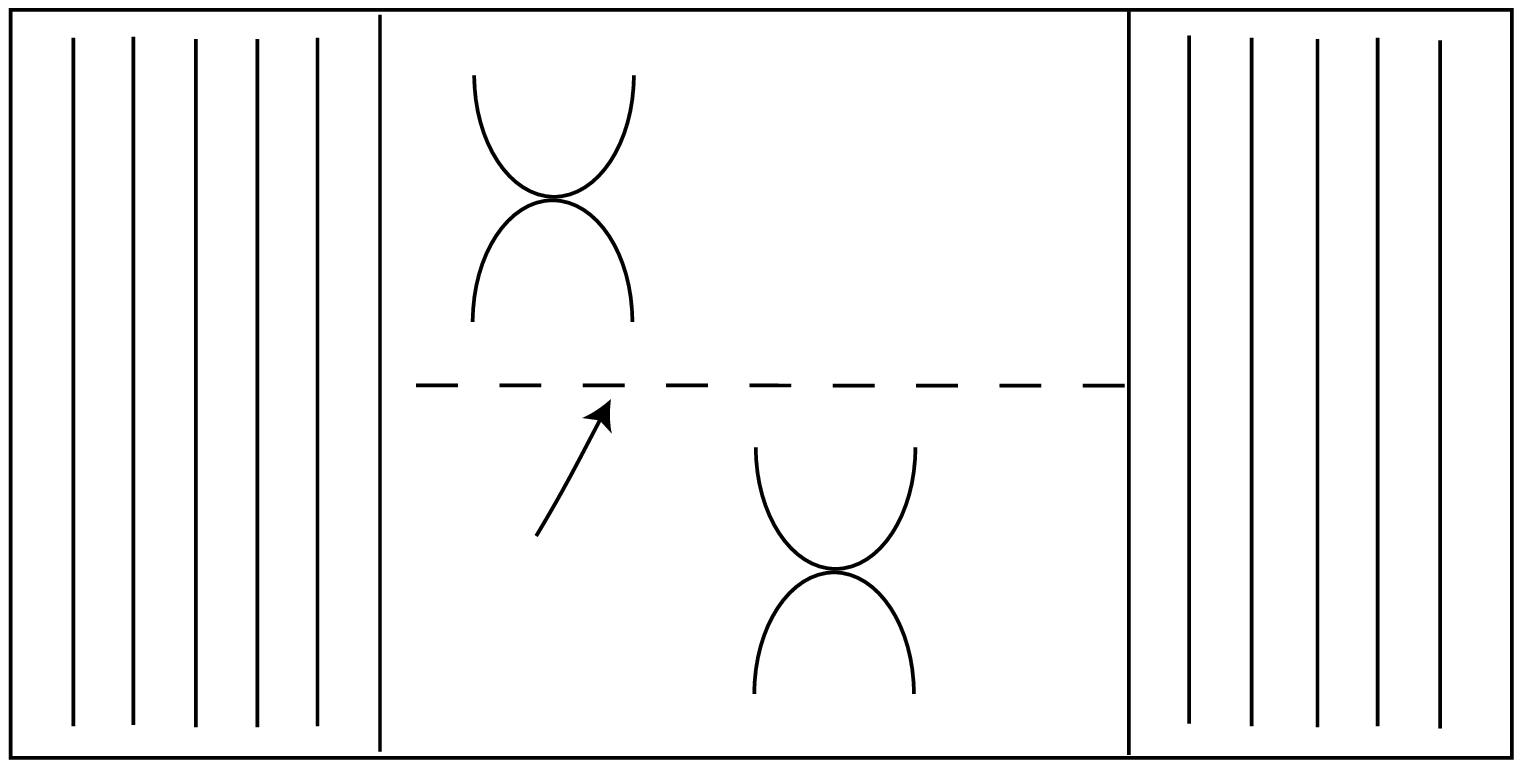}}}
 \caption{\label{3condfig}Summary of conditions $1,2,3$}
\end{figure}

Later, a comprehensive analysis of $1$--connected $h$--cobordisms
extended Akbulut's result (see \cite{CFHS,Ma,K1}). The following
picture of the general $1$--connected $5$--dimensional $h$--cobordism
$(W;P,Q)$ emerges. First, handles of indices $0,1,4$ and $5$ are
cancelled.  Let $L\subset W$ be the middle level between the
$2$--handles and $3$--handles. In $L$ lie the ascending $2$--spheres
$A$ of the $2$--handles
and the descending $2$--spheres $D$ of the $3$--handles.
It is possible to engulf
$A\cup D$ in a $4$--manifold $N\subset L$, where $N$ is homotopy
equivalent to a wedge of $2$--spheres.
The gradient lines through
$N$ define a sub-$h$--cobordism $(X;M,M')\subset (W;P,Q)$.
Combinatorially, $X$ is obtained from $N\times[-1,1]$ by attaching
$3$--handles to $D \times \{1\}$
and (up-side down) $3$--handles to $A \times \{-1\}$, so
$M=N/A$ and $M'=N/D$, where $/$ represents ``surgery".
($M$ and
$M'$ no longer denote the Mazur manifold, being instead
generalizations.)
The gradient lines constitue a product structure
on the complementary $h$--cobordism $(\overline{W\setminus
X};\overline{P\setminus M},\overline{Q\setminus M'})$, so
$(X;M,M')$ is the ``interesting" part.
By choosing $N$ with care,
the following conditions can be achieved:
\begin{enumerate}
\item $X\cong B^5$, so $M\cup \overline{M}'\cong S^4$.
\item The doubles $M\cup \overline{M}$ and $M'\cup \overline{M}'$
are both diffeomorphic to $S^4$.
\item There is a diffeomorphism $\tilde{\theta}\co M\rightarrow M'$,
so that $\tilde{\theta}|_{\bd M}$ composed with the gradient
flow identification $\bd M'\cong \bd M$ is an involution
$\theta\co \bd M \rightarrow \bd M$ (so $M\cup
\overline{M}'=M\cup_{\theta}\overline{M}$).
\end{enumerate}
Some of this information is summarized in \figref{3condfig}.

\begin{figure}[ht!]\anchor{N0pic}
\cl{
\small
\SetLabels
\E(.62*.9) $+$\\
\E(.62*.65) $-$\\
\E(.62*.35) $+$\\
\E(0*.5) $A$\\
\E(1*.5) $D$\\
\endSetLabels
\AffixLabels{
 \includegraphics[width=5cm, height=4cm]{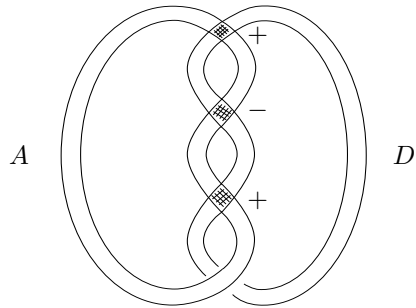}}}
 \caption{\label{N0pic}A schematic picture 
 of the neigbourhood $N_0$ of $A\cup D$}
\end{figure}

\begin{figure}[ht!]\anchor{N0kirby}
\cl{
\small
\SetLabels
\E(.8*.97) $s$\\
\E(.23*.02) $t$\\
\E(.04*.7) $0$\\
\E(.96*.7) $0$\\
\E(.15*.87) $A$\\
\E(.9*.85) $D$\\
\endSetLabels
\AffixLabels{
 \includegraphics[width=4.5cm]{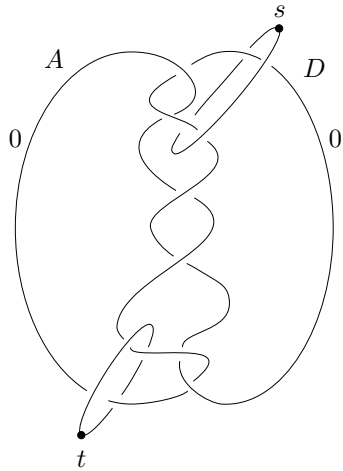}}}
 \caption{\label{N0kirby}A picture of $N_0$ in Kirby calculus notation}
\end{figure}

\begin{figure}[ht!]\anchor{+dotsfig}
\cl{
\small
\SetLabels
\E(.62*.97) $s$\\
\E(.16*.02) $t$\\
\E(.04*.78) $a$\\
\E(.73*.7) $0$\\
\E(.67*.85) $D$\\
\E(.97*.87) $u_1$\\
\E(.97*.67) $u_n$\\
\endSetLabels
\AffixLabels{
 \includegraphics[width=6cm]{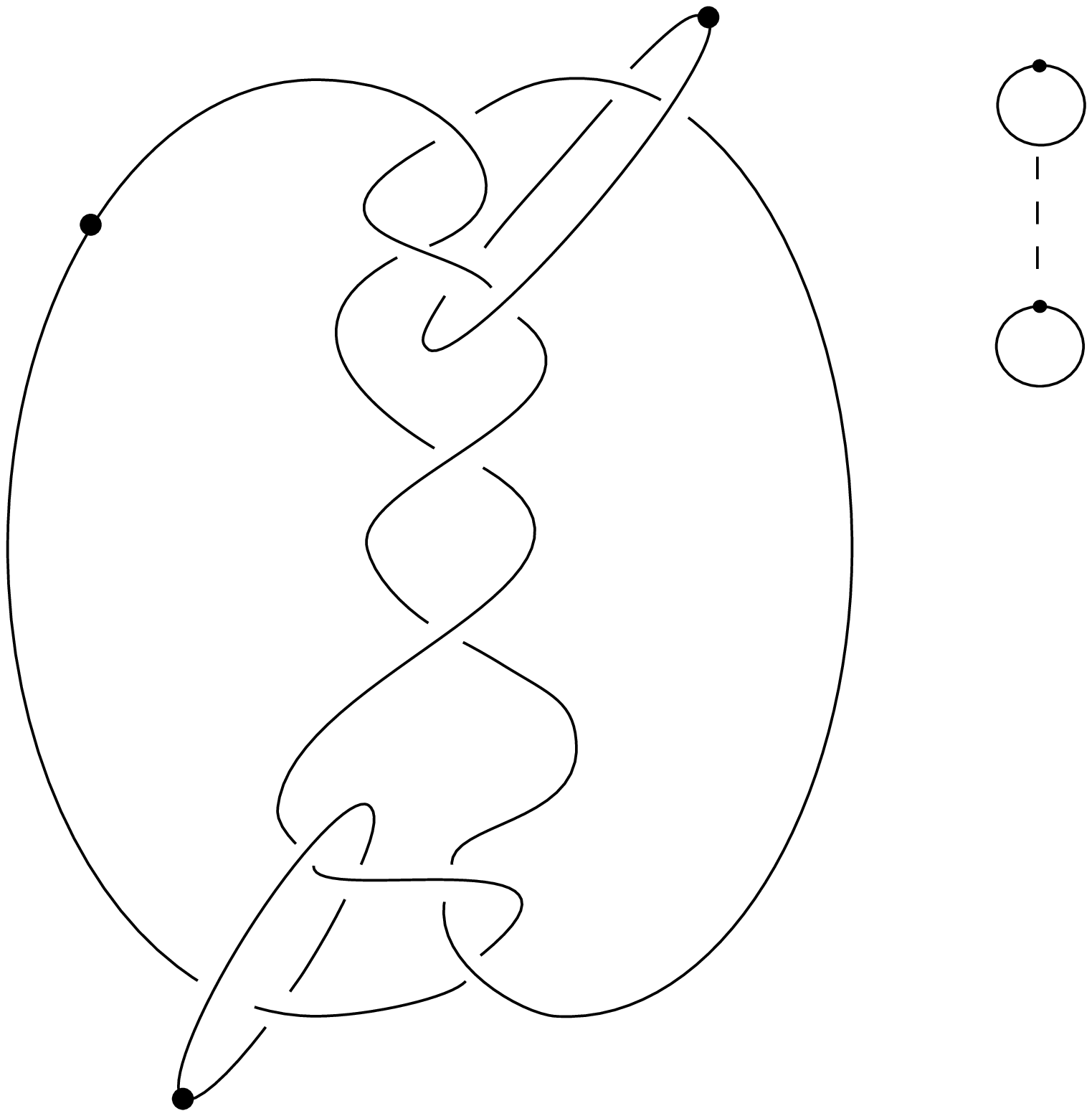}}}
 \caption{\label{+dotsfig}A picture of $M$ in Kirby calculus notation.
 Surgery on $A$ has converted one of the $2$--handles into a $1$--handle}
\end{figure}

We will explain point 2 above as a warm up to the non simply connected
case, Theorem~\ref{scobthm}.
We recall that $N$ is built from a neighbourhood
$N_0=\mathit{n}(A\cup D)$ of the ascending and descending spheres
arranged (with additional intersection points if necessary) so
that $\pi_1(L\setminus N_0)$ is trivial. A model picture of $N_0$ is
shown schematically in \figref{N0pic}. A more detailed
representation, using the Kirby calculus notation (see
\cite{Ak,K}) is given in \figref{N0kirby}. More complicated
configurations of $A\cup D$ require no new
ideas, only more notation, so we will treat the model case.
Fix a 
4--dimensional
handle decomposition $\MH$ (without $0$--handles) of
$(L\setminus N_0,\bd N_0)$. 
Form $N_1:=N_0\cup \{\mbox{all 1--handles of }\MH\}$ 
then 
$N=N_1\cup(\mbox{certain 2--handles})$. 
The $2$--handles are made by
stabilizing $\MH$ with $(2,3)$--handle pairs and passing the new
$2$--handles over suitable combinations of $2$--handles of $\MH$.
Since $\pi_1(L\setminus N_1)$ is trivial, we have complete freedom
in choosing the relations that these new $2$--handles determine
in $\pi_1(\bd(L \setminus N_1)) = \pi_1(\bd N_1)$. 
Since $\pi_1(\bd N_1) \to \pi_1(N_1)$ is an epimorphism we also
have complete freedom to add 2--handle relations to the presentation
of $\pi_1(N_1)$.
We will describe the relations that we introduce shortly.
Let $M_0 := N_0 / A$, $M_1 := N_1 / A$, $M := N / A$, and similarly
$M'_0 := N_0 / D$, etc.
Now
$M_0,M_1$ and $M$ result from surgery on $A$ and we represent this
diagrammatically by replacing the the $0$--framing at $A$ with a
dot, which indicates that the surgery has converted the $2$--handle
into a $1$--handle (see \figref{+dotsfig}). We proceed
similarly for $M_0',M_{1}'$ and $M'$.

To obtain $M \cup \overline{M}\cong S^{4}$, it is sufficient to show that
$M\times I\cong B^5$, the $5$--ball. But $M$ has just been
described as a $1,2$--handle body, so $M\times I$ is also a
$1,2$--handle body and in dimensions $d\ge 5$, only the group
theoretic presentation $\rho$ is relevant in determining if the
handle body is a ball:  $M\times I$ is a ball, $B^5$, if $\rho$ is
``deformable" to, or ``Andrews--Curtis related" to, the empty
presentation. The presentation $\rho$ that we may read off from
\figref{+dotsfig} has generators $a,s,t,u_1,\ldots,u_n$ and so
far only one relation: $t^{-1}ata^{-1}s^{-1}as$. This length $7$
relation may not look like a promising start for an a standard
presentation of the trivial group, but we can begin by choosing
two new $2$--handles representing $t$ and $s$, which collapse the
relation to $a$. From here, simply choose $2$--handles representing
length $1$ relations $u_1,\ldots,u_n$. In this way, $N$ is built
from $N_1$ so that both $M \times I\cong B^5$ and $M' \times I\cong
B^5$. For the construction of the involution $\theta$, see
lemma~$2$ in \cite{Ma}. The statement that $X\cong B^5$ in point
$1$ above can be extracted from our proof of theorem~\ref{scobthm}
in the case $k=0$.

We can now prove the following:
\begin{thm}
\label{homthm} UTQFTs cannot distinguish $1$--connected smooth
$4$--manifolds which are homotopy equivalent. In fact, even a less
rigid ``theory", where Atiyah's gluing axiom is only enforced
along homology $3$--spheres, will likewise be unable to distinguish
homotopy equivalent $4$--manifolds.
\end{thm}
\proof It is well known that smooth homotopy equivalent
$4$--manifolds $P$ and $Q$ are smoothly $h$--cobordant. In the preceding
notation write $P=C\cup_S M$ and $Q=C\cup_S M'$, where $S=\bd
M$. Let $Z$ be a UTQFT and let $V(S)$ be the Hilbert space
assigned to $S$ by $Z$, so $Z(C),Z(M),Z(M')\in V(S)$. Now $\langle
M-M', M-M'\rangle = S^4-S^4-S^4+S^4=0\in\M^4$, so $\langle
Z(M-M'), Z(M-M')\rangle =0\in V(S)$ and by unitarity,
$Z(M-M')=0\in V_S$. Hence $Z(M)=Z(M')$. Finally, $Z(P)=\langle
Z(C),Z(M)\rangle =\langle Z(C),Z(M')\rangle = Z(Q)$. \endproof

We may extend this result to the non-simply connected setting:
\begin{thm}
\label{scobthm} UTQFTs\,cannot\,distinguish\,smoothly $s$-cobordant\,%
$4$-manifolds.
\end{thm}
{\bf Remark}\qua The proof is a mild generalization of the
preceding middle level analysis. It should be noted that the
$3$--manifolds whose Hilbert spaces we must now consider are no
longer homology $3$--spheres, but instead admit maps to
$\underset{k}{\vee} S^{1}$ and have vanishing $H_1(\;\cdot\; ;\mathbb{Z}
[\pi_1(\underset{k}{\vee} S^1)])$, ie, the corresponding covers have perfect
fundamental groups.

\proof[Proof of Theorem~\ref{scobthm}] We assume
our manifolds are oriented. Again, let $(W;P,Q)$ be the
$s$--cobordism with $0,1,4$ and $5$--handles eliminated. Let $L$,
again, be the middle level. Find ascending spheres $A = \cup A_i$ and
descending spheres $D = \cup D_j$ for a handle structure exhibiting
triviality in the Whitehead group, and so that $\pi_1(L\setminus
(A\cup D))\rightarrow \pi_1(L)$ is an isomorphism.
This may require 
stabilizations and
``finger moves" between the original $A$ and $D$.

The use of finger moves --- the inverse Whitney's famous trick \cite{Wh} ---
to improve fundamental groups goes back to Casson's constructions \cite{C}.
The particular incarnation used here is explained in Section 7.1.D of \cite{FQ},
the first half of page 105 being the crux.

Our goal is to engulf $A\cup D$ in $N\subset L$ so that there is a
corresponding sub-$s$--cobordism $(X;M,M')$ so that the complement
has the gradient product structure and
\begin{enumerate}
\item $X\cong \underset{k}{\natural} (S^1\times B^4)$, so
$M\cap \overline{M}'\cong \underset{k}{\sharp} (S^1\times S^3)$.
\item The doubles $M\cup \overline{M}$ and $M'\cup \overline{M}'$
are both diffeomorphic to $\underset{k}{\sharp} (S^1\times S^3)$.
\end{enumerate}
By choosing $\Delta$ symmetrically below we may also arrange that
$M\cong M'$ and that the corresponding $\theta\co \bd
M\rightarrow\bd M$ is an involution, but these statements are
irrelevant to the conclusions on UTQFTs
so we do not expatiate.

\begin{figure}[ht!]\anchor{deltafig}
\cl{
\small
\SetLabels
\E(.2*.6) $\Delta$\\
\E(.1*.95) $A$\\
\E(.9*.96) $D$\\
\endSetLabels
\AffixLabels{
 \includegraphics[width=6cm, height=4.6cm]{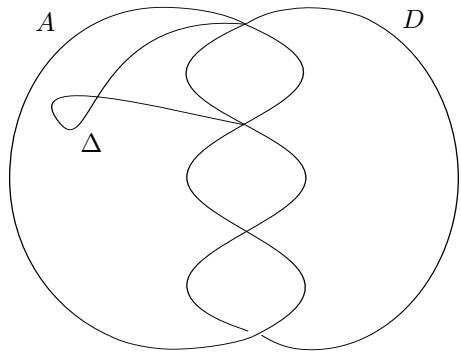}}}
 \caption{\label{deltafig}A picture of $A\cup D\cup \Delta$}
\end{figure}

Since we are trying to build an $N$ with $\pi_1(N)$ a large free
group, rather than a trivial group, we will not have to be as
careful in enlarging $A\cup D$ to $N$ as before. Let $\Delta$ be a
union of immersed Whitney disks for $A\cap D$ which pair ``excess"
double points (see \figref{deltafig}).
Set $N={\mathit n}(A\cup D\cup \Delta)$. Now $X$ is a thickening
of $D_{A}^{3}\cup D_{D}^{3}\cup \Delta$, where $D_{A}^{3}$ and
$D_{D}^{3}$ are the $3$--disks descending from $A$ and ascending
from $D$, so $X$ manifestly collapses to a wedge of circles. Thus,
$X\cong \underset{k}{\natural}(S^{1}\times B^{4})$, since $W$ is
orientable.

\begin{figure}[ht!]\anchor{Xfig}
\cl{
\small
\SetLabels
\E(.1*.82) $a$\\
\E(.9*.8) $D$\\
\E(.98*.6) $0$\\
\E(.7*.97) $x$\\
\E(.83*.0) $y$\\
\endSetLabels
\AffixLabels{
 \includegraphics[width=6.5cm, height=6.5cm]{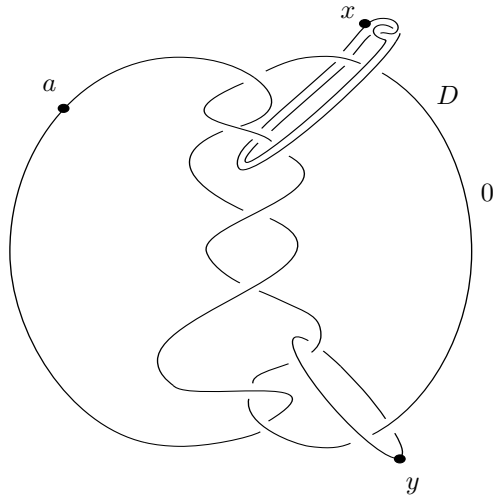}}}
 \caption{\label{Xfig}Model diagram for $M$}
\end{figure}

The final point is that $D(M)\cong
D(M')=\underset{k}{\natural}S^1\times S^3$. As before, this
reduces to seeing a $1$--handle presentation for $M$ ($M'$) and
verifying that that presentation deforms to the standard
presentation of the $k$--generator free group. The model diagram
for $M$ is given in \figref{Xfig}.

The reader should compare this with \figref{N0pic}.
The first thing we notice is that attaching a 3--handle 
(which effects a surgery on $A$)
has been described by removing a zero label and adding a dot, indicating
a 1--handle.
Locally we are simply creating a bordism from $S^2 \times D^2$ to 
$D^3 \times S^1$ (rel $S^2 \times S^1$) by attaching a 3--handle
$(D^3 \times D^2, S^2 \times \frac{1}{2} D^2)$ to
$S^2 \times \frac{1}{2} D^2 \subset S^2 \times D^2$.
The bottom of the bordism $S^2 \times D^2$ is represented by a trivial
2--handle, the top by a 1--handle.

The corresponding model diagram for $M'$ is obtained by replacing
the 1--handle $a$ by a 0--framed 2--handle $A$ and
the 0--framed 2--handle $D$ by a 1--handle $d$.
The Whitehead double (curve $x$ in the figure)
arises from the attachment of $\Delta$;
there will be one such $x$ for each double point on $\Delta$,
whereas only a single $y$ curve is present for each pair $(A_i, D_i)$.
The boundary of the
$2$--handle core $D$ reads ``$a$", so the presentation is a standard
one for the free group: $\{a,x,y: a\}$. To see this presentation,
find disjoint surfaces $\bar{a}$, $\bar{x}$ and $\bar{y}$ bounding
$a,x,y$ such that $D\cap(\bar{a}\cup\bar{x}\cup\bar{y})$ is just
one point lying on $\bar{a}$. 
As $D$ crosses $\bar{a}$, $\bar{x}$ and $\bar{y}$ it reads its relation in
the free group generated by $a$, $x$ and $y$;
the result is simply a smaller free group.
The general case,
involving additional stabilizations and more double points, is similar.

\section{Problems}

{\bf Problem 1a}\qua Given $x\in \M^{d}_{S}$ with $\langle
x,x \rangle \neq 0$, is it possible to construct a UTQFT which
assigns to $x$ a nonzero vector $\tilde{x}\in V_{S}$?
\medskip

{\bf Problem 1b}\qua Similarly, given $x\in \M^{d}_{S}$ such that
all $d$--dimensional UTQFTs assign the zero vector to $x$, does it
follow that $\langle x,x \rangle=0$? \medskip

{\bf{Problem 2}}\qua Is there a $3$--manifold $S$ and a
nonzero vector $x\in \M^4_S$ such that $\langle x,y \rangle =0$
for all $y \in \M^4_S$? 
Such an $x$ would be called ``singular".

A positive answer to this question would give an example
of a combination of $4$--manifolds that will be sent to $0$ in any
TQFT with a non singular quadratic form, not just in unitary TQFTs.
(In a UTQFT $Z$, $\langle Z(x),Z(x) \rangle = 0$ implies that
$\langle Z(x),Z(y) \rangle=0$ for all $y$.)
\medskip

{\bf{Problem 3}}\qua Analyze the $3$--dimensional
pairings.\medskip

{\bf Problem 4}\qua Analyze the pairings in dimensions
$d\ge 5$ .\medskip

{\bf Problem 5}\qua In dimension $d=4$ characterize the
zero locus of the pairing in more detail. For instance, are there
any elements in the zero locus of the pairing that have an odd
number of manifolds with non-zero coefficients?

Even numbers of nonzero coefficients may be obtained using
manifolds which are disjoint unions of the examples given already
\medskip

{\bf Problem 6}\qua Consider the same problem in the
piecewise linear and topological categories.
\medskip

{\bf{Problem 7}}\qua Consider coefficients other than
$\C$.

This is almost certainly of interest,
because there exist classes of TQFTs whose invariants take values
in rings other than $\C$. Clearly the pairings will never be
positive in any dimension for coefficient rings with elements $x$
with satisfy $x\bar{x}=0$, but nevertheless, even for such rings,
a characterisation of the nullity may be interesting. For example,
in the ring $(\Z/7\Z)[\omega]$ with $\omega=e^{2\pi i/7}$ and the
involution given by extension of $\bar{1}=1$ and
$\bar{\omega}=\omega^{-1}$, there are elements $a$ for which
$a\bar{a}=0$, but using the Milnor sphere of order $7$, one can
construct more interesting examples $x \in
\M_{S^6}[(\Z/7\Z)[\omega]]$ which have $\langle x,x\rangle=0\in
\M[(\Z/7\Z)[\omega]]$, such as $x=\sum_{i=0}^6\omega^i (B^6,\theta_i)$,
where $\theta$ is the ``clutching map" for the Milnor sphere.

\end{document}